\documentclass[a4paper,11pt]{amsart}
\usepackage{graphicx}
\usepackage{mathptmx}
\usepackage{amsmath}
\usepackage{amssymb}
\usepackage{enumitem}
\usepackage{xcolor}
\usepackage{pgfplots}

\newmuskip\pFqmuskip

\newcommand*\pFq[6][8]{%
  \begingroup 
  \pFqmuskip=#1mu\relax
  \mathcode`=\string"8000
  \begingroup\lccode`\~=`\,
  \lowercase{\endgroup\let~}\pFqcomma
  F^{#2}_{#3}{\left(\genfrac..{0pt}{}{#4}{#5}\bigg|#6\right)}%
  \endgroup
}
\newcommand{\pFqcomma}{\mskip\pFqmuskip}

\newcommand{\ap}{a^{\dagger}}
\newcommand{\ba}{aa^{\dagger}}
\newcommand{\fa}{a^{\dagger}a}

\newtheorem{theorem}{Theorem}[section]

\begin{document}

\title[Degenerate generalized Stirling operators of the first kind]{Degenerate generalized Stirling operators of the first kind arising from generalized Heisenberg algebra}

\author{Taekyun  Kim}
\address{Department of Mathematics, Kwangwoon University, Seoul 139-701, Republic of Korea}
\email{tkkim@kw.ac.kr}
\author{Dae San  Kim}
\address{Department of Mathematics, Sogang University, Seoul 121-742, Republic of Korea}
\email{dskim@sogang.ac.kr}

\subjclass[2010]{11B73; 11B83}
\keywords{degenerate generalized Stirling operators of the first kind; degenerate generalized $r$-Stirling operators of the first kind; generalized Heisenberg algebra}

\begin{abstract}
This paper investigates the degenerate generalized Stirling operators of the first kind, $S_{f,\lambda}(n,k,H)$, bridging a gap in the operational calculus of the generalized Heisenberg algebra $\text{GHA}_f$ unified with degenerate calculus. As they are the inverse of the degenerate generalized Stirling operators of the second kind, these operators express the monomial operator products $(a^\dagger)^n a^n$ in terms of the degenerate factorial operators $(a^\dagger a)_{k,\lambda}$. We derive key structural and combinatorial properties for these operators, including an explicit product factorization, a fundamental recurrence relation, and an operational shifting identity. Furthermore, we establish the orthogonality relations between the degenerate generalized Stirling operators of the first and second kinds, providing a complete combinatorial framework for functional quantum algebras.
\end{abstract}

\maketitle

\section{Introduction}
The study of combinatorics intersecting with quantum mechanics has a rich history, particularly in the operational calculus of the Heisenberg-Weyl algebra. Central to this intersection is the problem of normal ordering, which involves rewriting functions of the creation operator $\ap$ and annihilation operator $a$ into a standard form where all creation operators precede all annihilation operators. Traditionally, the coefficients arising from the normal ordering of powers of the number operator $N = \fa$ are given by the classic Stirling numbers of the second kind. \par
In recent years, significant attention has been paid to deformations and generalizations of these algebraic structures to better model complex physical systems. Two prominent avenues of generalization include the degenerate calculus, which introduces a parameter $\lambda$ that alters the underlying difference operators, and the generalized Heisenberg algebra $\text{GHA}_{f}$. The latter replaces standard commutation relations with functional forms dependent on a Hermitian operator $H$ and an analytic function $f(H)$, offering a highly flexible framework that encompasses various deformed quantum systems. While normal ordering and Stirling numbers of the second kind have been extensively explored within these generalized settings, the corresponding inverse relations -- specifically, expressing powers of monomial operators in terms of factorial or degenerate operational variants -- remain less developed. The primary objective of this paper is to bridge this gap by systematically investigating the degenerate generalized Stirling operators of the first kind. By establishing their fundamental recurrence relations, explicit formulas, and algebraic properties, we provide a complete combinatorial framework that unifies degenerate calculus with functional quantum algebras. \par
The remainder of this paper is organized as follows. In Section 2, we establish the necessary mathematical preliminaries. We recall the classic Stirling numbers of the first kind and the degenerate Stirling numbers of both kinds, followed by a reminder of the degenerate boson normal ordering expansion of $(\fa)_{n,\lambda}$ within the standard Heisenberg algebra. We then review the generalized Heisenberg algebra $\text{GHA}_{f}$, which is generated by the set $\{a,\ap,H\}$ and governed by the commutation relations:
\begin{equation*}
aH=f(H)a,\quad H\ap=\ap f(H),\quad \big[a,\ap\big]=a \ap-\ap a=f(H)-H,
\end{equation*}
where $f(H)$ is an analytic function of $H$. To facilitate our derivations, we introduce the quantum operator functions $[k]_{f}(H)=f^{k}(H)-H$ for any integer $k$, obtaining several structural relations involving these operators alongside powers of the creation and annihilation operators. After reviewing the operational actions of $a$, $\ap$, and $N_{f}=\fa$ on the number states $|m\rangle$, we recall the generalized Stirling operators of the second kind and their degenerate counterparts, defined via the normal ordering expansion for $n \ge 0$:
\begin{equation}
\big(N_{f}\big)_{n,\lambda}=\big(\fa\big)_{n,\lambda}=\sum_{k=0}^{n}{n \brace k}_{f,\lambda}(H)\big(\ap\big)^{k}a^{k}. \label{0}
\end{equation}
In Section 3, we turn our attention to the primary focus of this work: the degenerate generalized Stirling operators of the first kind, $S_{f,\lambda}(n,k,H)$. These operators are defined via the inverse relation of $\eqref{0}$ for $n \ge 0$:
\begin{equation*}
\big(\ap\big)^{n}a^{n} = \sum_{k=0}^{n}S_{f,\lambda}(n,k,H)\big(\fa\big)_{k,\lambda}.
\end{equation*}
We present several core structural properties of these operators: Theorem 3.1 establishes the explicit factorization product for any positive integer $k$: $(\ap)^{k}a^{k}=\prod_{l=0}^{k-1}\big(\fa-(l)_{f}(H)\big)$.
Theorem 3.2 derives a fundamental recurrence relation for any integers $n \ge k \ge 0$:\begin{equation*}
S_{f,\lambda}(n+1,k,H)=S_{f,\lambda}(n,k-1,H)+\big(k\lambda-(n)_{f}(H)\big)S_{f,\lambda}(n,k,H).
\end{equation*}
Theorem 3.3 demonstrates that the degenerate generalized Stirling operators of the first and second kind satisfy standard orthogonality relations. Theorem 3.4 yields a general operational identity for any nonnegative integers $n$ and $r$:
\begin{equation*}
(a^\dagger)^{r}(\fa)_{n}= \sum_{k=0}^{n}S_{f,\lambda}^{(r)}\big(n+r,k+r,f^{-r}(H)\big)\big(\fa\big)_{k,\lambda}\big(\ap\big)^{r}.
\end{equation*}
Finally, Section 4 provides a brief summary and concludes the paper.

\section{Preliminaries}
Let $\lambda$ be a fixed nonzero real number. The degenerate falling factorial sequence is given by
\begin{equation*}
(x)_{0,\lambda}=1,\quad (x)_{n,\lambda}=x(x-\lambda)(x-2\lambda)\cdots\big(x-(n-1)\lambda\big),\ (n\ge 1),\ (\mathrm{see}\ [7-12]).
\end{equation*}
The classical Stirling numbers of the first kind are defined by
\begin{equation*}
(x)_{n}=\sum_{k=0}^{n}S_{1}(n,k)x^{k},\quad (n\ge 0),\quad (\mathrm{see}\ [7-12]),
\end{equation*}
where
\begin{equation*}
(x)_{0}=1,\quad (x)_{n}=x(x-1)(x-2)\cdots (x-n+1),\ (n\ge 1).
\end{equation*}
In [11], the degenerate Stirling numbers of the first kind are given by
\begin{equation}
(x)_{n}=\sum_{k=0}^{n}S_{1,\lambda}(n,k)(x)_{k,\lambda},\quad (n\ge 0),\quad (\mathrm{see}\ [11]).\label{1}
\end{equation}
As the inversion formula of \eqref{1}, the degenerate Stirling numbers of the second kind are defined by
\begin{equation}
(x)_{n,\lambda}=\sum_{k=0}^{n}{n \brace k}_{\lambda}(x)_{n},\quad (n\ge 0),\quad (\mathrm{see}\ [7,10,11]).\label{2}
\end{equation} \par
The degenerate boson normal ordering is given by the creation $\ap$ and annihilation $a$ operators as follows:
\begin{equation}
\big(\fa\big)_{n,\lambda}=\sum_{k=0}^{n}{n \brace k}_{\lambda}\big(\ap\big)^{k}a^{k},\quad (n\ge 0),\quad (\mathrm{see}\ [8,10,11]). \label{3}	
\end{equation}
The Heisenberg algebra consists of the algebra generated by $\{a,\ap\}$ satisfying the commutation relations:
\begin{equation*}
aH=(1+H)a,\quad H\ap=\ap\big(1+H),\quad\textrm{and}\quad [a,\ap]=\ba-\fa=1,
\end{equation*}
where $H=\fa$ is the Hamiltonian. \par
Recently, Ouahhabi and Tahri introduced the generalized Heisenberg algebra, denoted by $\text{GHA}_{f}$, and quantum operator functions. The generalized Heisenberg algebra consists of the algebra generated by $\{a,\ap,H\}$ such that the following commutation relations hold:
\begin{equation}
aH=f(H)a,\quad H\ap=\ap f(H),\label{4}
\end{equation}
and
\begin{equation}
[a,\ap]=\ba-\fa=f(H)-H,\label{5}
\end{equation}
where $H$ is the self-conjugate Hamiltonian of the system and $f(H)$ is an analytic function of $H$ (see [14]). For any positive integer $k$, the power compositions of $f$ are defined by
\begin{equation*}
f^{k}=\underbrace{f\circ f\circ \cdots \circ f}_{k-\mathrm{times}},\quad f^{-k}=\underbrace{f^{-1}\circ f^{-1}\circ \cdots \circ f^{-1}}_{k-\mathrm{times}},
\end{equation*}
where $f^{-1}$ is the compositional inverse function of $f$ and $f^{0}$ denotes the identity function. \\
Let $k\in\mathbb{Z}$. Then the quantum operator functions of the first kind are defined by
\begin{equation}
[k]_{f}(H)=f^{k}(H)-H,\quad (\mathrm{see}\ [14]). \label{6}
\end{equation}
Also, the quantum operator functions of the second kind are given by
\begin{equation}
(k)_{f}[H]=[k]_{f}\big(f^{-k}(H)\big)=H-f^{-k}(H).\label{7}
\end{equation}
It is known that
\begin{equation}
a^{m}[n]_{f}\big(f^{j}(H)\big)=[n]_{f}\big(f^{m+j}(H)\big)a^{m},\ [n]_{f}\big(f^{j}(H)\big)(\ap)^{m}=(\ap)^{m}[n]_{f}\big(f^{j+m}(H)\big), \label{8}
\end{equation}
and
\begin{equation}
\big[a,(\ap)^{k}\big]=(\ap)^{k-1}[k]_{f}(H),\quad [a^{k},\ap]=[k]_{f}(H)a^{k-1},\quad (\mathrm{see}\ [9,14]).\label{9}
\end{equation}
In particular, from \eqref{8}, we see that
\begin{equation}
[k]_{f}(H)\fa=\fa[k]_{f}(H). \label{10}
\end{equation} \par
Let
\begin{equation}
H\big|m\big\rangle=\alpha_{m}\big|m\big\rangle,\label{11}
\end{equation}
where $\alpha_{m}$ is the dimensionless energy spectrum of the quantum system. \\
Then, by \eqref{11}, we get
\begin{equation*}
\alpha_{m}=f(\alpha_{m-1})=f^{2}(\alpha_{m-2})=\cdots=f^{m}(\alpha_{0}),\quad (\mathrm{see}\ [9,14]).
\end{equation*}
Let $N_{f}=\fa$. Then we have
\begin{equation}
N_{f}\big|m\big\rangle=\fa\big|m\big\rangle=[m]_{f}(\alpha_{0})\big|m\big\rangle.\label{12}
\end{equation}
From \eqref{12}, the actions of $a$ and $\ap$ operators on the Fock space representation of the generalized Heisenberg algebra are given by
\begin{equation}
a|m\rangle=\sqrt{[m]_{f}(\alpha_{0})}\Big|m-1\Big\rangle,\quad a|0\rangle=0, \label{13}	
\end{equation}
and
\begin{equation}
\ap|m\rangle=\sqrt{[m+1]_{f}(\alpha)_{0}}\Big|m+1\Big\rangle,\quad (m\ge 0),\quad (\mathrm{see}\ [9,14]). \label{14}
\end{equation} \par
Recently, Ouahhabi and Tahri introduced the generalized Stirling operators of the second kind as
\begin{equation*}
N_{f}^{n}=\big(\fa\big)^{n}=\sum_{k=0}^{n}{n \brace k}_{f}(H)\big(\ap\big)^{k}a^{k},\quad (n\ge 0),\quad (\mathrm{see}\ [9,14]). 	
\end{equation*}
Note that
\begin{equation*}
N_{1+H}^{n}=\big(\fa\big)^{n}=\sum_{k=0}^{n}{n \brace k}_{1+H}(H)\big(\ap\big)^{k}a^{k}=\sum_{k=0}^{n}{n \brace k}\big(\ap\big)^{k}a^{k},
\end{equation*}
where ${n \brace k}$ are the Stirling numbers of the second.
In [9], the degenerate generalized Stirling operators of the second kind are defined by (see \eqref{2}, \eqref{3})
\begin{equation}
\big(N_{f}\big)_{n,\lambda}=\big(\fa\big)_{n,\lambda}=\sum_{k=0}^{n}{n \brace k}_{f,\lambda}(H)\big(\ap\big)^{k}a^{k},\quad (n\ge 0),\quad (\mathrm{see}\ [14]). \label{15}
\end{equation}
Note that
\begin{equation*}
{n \brace k}_{H+1,\lambda}(H)={n \brace k}_{\lambda},\quad (n,k\ge 0),\quad (\mathrm{see}\ [9]).
\end{equation*}
General references of this paper include [1-3,5,6,13,15-18].

\section{Degenerate generalized Stirling operators of the first kind arising from generalized Heisenberg algebra}
In this section, we study the degenerate generalized Stirling operators of the first kind arising from the inverse relation of \eqref{15}.\\
From \eqref{7}-\eqref{10}, we note that
\begin{align}
\big(\ap\big)^{k}a&=a(\ap)^{k}-(\ap)^{k-1}[k]_{f}(H) \label{16}\\
&=a\big(\ap\big)^{k}-[k]_{f}\big(f^{-k+1}(H)\big)\big(\ap\big)^{k-1} \nonumber\\
&=\Big(\ba-[k]_{f}\big(f^{-k+1}(H)\big)\Big)\big(\ap\big)^{k-1}	\nonumber\\
&=\Big(\fa+f(H)-H-f(H)+f^{-k+1}(H)\Big)\big(\ap\big)^{k-1}\nonumber\\
&=\big(\fa-(k-1)_{f}(H)\big)\big(\ap\big)^{k-1},\nonumber
\end{align}
and, by \eqref{16}, we get
\begin{equation*}
\begin{aligned}
\big(\ap\big)^{k}a^{2}&=\Big(\fa-(k-1)_{f}(H)\Big)\big(\ap\big)^{k-1}a\\
&=\Big(\fa-(k-1)_{f}(H)\Big)\Big(\fa-(k-2)_{f}(H)\Big)\big(\ap\big)^{k-2}.
\end{aligned}
\end{equation*}
Continuing this process, we have  the following theorem. Note here that, because of \eqref{10}, the ordering of the factors in \eqref{17} doesn't matter.
\begin{theorem}
For any positive integer $k$, we have
\begin{align}
(\ap)^{k}a^{k}&=\big(\fa-(k-1)_{f}(H)\big)\big(\fa-(k-2)_{f}(H)\big)\cdots\big(\fa-(1)_{f}H\big)\fa \label{17}\\
&=\prod_{l=0}^{k-1}\big(\fa-(l)_{f}(H)\big). \nonumber
\end{align}
\end{theorem}
As the inverse relation of \eqref{15}, the degenerate generalized Stirling operators of the first kind are defined in [4] by (see \eqref{1})
\begin{equation}
\big(\ap\big)^{n}a^{n}	=\sum_{k=0}^{n}S_{f,\lambda}(n,k,H)\big(\fa\big)_{k,\lambda},\quad (n\ge 0). \label{18}
\end{equation}
From Theorem 2.1 and \eqref{18}, we note that
\begin{equation*}
S_{f,\lambda}(n,k,H)=0\ \textrm{if $k>n$,}\quad S_{f,\lambda}(n,0,H)=\delta_{n,0}, \quad S_{f,\lambda}(n,n,H)=1,\quad (n\ge 0).
\end{equation*}
where $\delta_{n,0}$ is Kronecker delta. \\
From \eqref{4}, we note that
\begin{equation} \label{19}
\begin{aligned}
a^{m}S_{f,\lambda}\big(n,k,f^{j}(H)\big)&=S_{f,\lambda}\big(n,k,f^{j+m}(H)\big)a^{m},\\
S_{f,\lambda}\big(n,k,f^{j}(H)\big)\big(\ap\big)^{m}&=\big(\ap\big)^{m}S_{f,\lambda}\big(n,k,f^{j+m}(H)\big).
\end{aligned}	
\end{equation}
In particular, from \eqref{19} we observe that
\begin{equation}
N_{f}S_{f,\lambda}(n,k,H)=S_{f,\lambda}(n,k,H)N_{f}. \label{19-1}
\end{equation} \par
From \eqref{18}, we have
\begin{equation}
\big(\ap\big)^{n+1}a^{n+1}=\sum_{k=0}^{n+1}S_{f,\lambda}(n+1,k.H)\big(\fa\big)_{k,\lambda}.\label{20}
\end{equation}
By \eqref{5} and \eqref{8}, we get
\begin{align}
&\big(\ap\big)^{n}a=\big(\ap\big)^{n-1}\fa\label{21}\\
&=\big(\ap\big)^{n-1}\big(\ba-[1]_{f}(H)\big)=\big(\ap\big)^{n-1}\ba-\big(\ap\big)^{n-1}[1]_{f}(H)\nonumber	\\
&=\big(\ap\big)^{n-1}\ba-[1]_{f}\Big(f^{-n+1}(H)\Big)\big(\ap\big)^{n-1}\nonumber\\
&=\Big[\big(\ap\big)^{n-2}\ba-[1]_{f}\big(f^{-n+2}(H)\big)\big(\ap\big)^{n-2}\Big]\ap-[1]_{f}\big(f^{-n+1}(H)\big)\big(\ap\big)^{n-1} \nonumber\\
&=\big(\ap\big)^{n-2}a\big(\ap\big)^{2}-[1]_{f}\big(f^{-n+2}(H)\big)\big(\ap\big)^{n-1}-[1]_{f}\big(f^{-n+1}(H)\big)\big(\ap\big)^{n-1} \nonumber\\
&=\cdots \nonumber\\
&=a\big(\ap\big)^{n}-\Big([1]_{f}(H)+[1]_{f}\big(f^{-1}(H)\big)+\cdots+[1]_{f}\big(f^{-n+1}(H)\big)\Big)\big(\ap\big)^{n-1}. \nonumber
\end{align}
From \eqref{6}-\eqref{8} and \eqref{21}, we note that
\begin{align}
&\big(\ap\big)^{n+1}a^{n+1}=\ap\big((\ap)^{n}a\big)a^{n}\label{22} \\	
&=\fa\big(\ap\big)^{n}a^{n}\nonumber\\
&\quad -\ap\big[[1]_{f}(H)+[1]_{f}\big(f^{-1}(H)\big)+\cdots+[1]_{f}\big(f^{-n+1}(H)\big)\big](\ap)^{n-1}a^{n}\nonumber\\
&=\fa\big(\ap\big)^{n}a^{n}\nonumber\\
&\quad -\big[[1]_{f}(f^{-1}(H))+[1]_{f}\big(f^{-2}(H)\big)+\cdots+[1]_{f}\big(f^{-n}(H)\big)\big](\ap)^{n}a^{n}\nonumber\\
&=\big(\fa\big)\big(\ap\big)^{n}a^{n}-(n)_{f}(H)\big(\ap\big)^{n}a^{n}.\nonumber
\end{align}
Thus, by \eqref{18}, \eqref{19-1} and \eqref{22}, we get
\begin{align}
&\big(\ap\big)^{n+1}a^{n+1}=N_{f}\sum_{k=0}^{n}S_{f,\lambda}(n,k,H)\big(\fa\big)_{k,\lambda}-(n)_{f}(H)\sum_{k=0}^{n}S_{f,\lambda}(n,k,H)\big(\fa\big)_{k,\lambda} \label{23} \\
&=\sum_{k=0}^{n}S_{f,\lambda}(n,k,n)\big(\fa-k\lambda+k\lambda\big)\big(\fa\big)_{k,\lambda}-(n)_{f}(H)\sum_{k=0}^{n}S_{f,\lambda}(n,k,H)\big(\fa\big)_{k,\lambda}\nonumber\\
&=\sum_{k=0}^{n}S_{f,\lambda}(n,k,H)\big(\fa\big)_{k+1,\lambda}+\sum_{k=0}^{n}\big(k\lambda-(n)_{f}(H)\big)S_{f,\lambda}(n,k,H)\big(\fa\big)_{k,\lambda}\nonumber\\
&=\sum_{k=0}^{n+1}\big\{S_{f,\lambda}(n,k-1,H)+\big(k\lambda-(n)_{f}(H)\big)S_{f,\lambda}(n,k,H)\big\}(\fa)_{k,\lambda}.\nonumber
\end{align}
Therefore, by \eqref{20} and \eqref{23}, we obtain the following recurrence relation for the degenerate generalized Stirling operators of the first kind.
\begin{theorem}
For any integers $n,k$ with $n \ge k \ge 0$, we have
\begin{equation*}
S_{f,\lambda}(n+1,k,H)=S_{f,\lambda}(n,k-1,H)+\big(k\lambda-(n)_{f}(H)\big)S_{f,\lambda}(n,k,H).
\end{equation*}
\end{theorem}
On the one hand, from \eqref{15} and \eqref{18}, we have
\begin{align}
\big(\ap\big)^{n}a^{n}&=\sum_{k=0}^{n}S_{f,\lambda}(n,k,H)\big(\fa\big)_{k,\lambda}=\sum_{k=0}^{n}S_{f,\lambda}(n,k,H)\sum_{j=0}^{k}{k \brace j}_{f,\lambda}(H)\big(\ap\big)^{j}a^{j}\label{24} \\
&=\sum_{j=0}^{n}\bigg(\sum_{k=j}^{n}S_{f,\lambda}(n,k,H){k \brace j}_{f,\lambda}(H)\bigg)\big(\ap\big)^{j}a^{n}.\nonumber
\end{align}
On the other hand, we also get
\begin{align}
\big(\fa\big)_{n,\lambda}&=\sum_{k=0}^{n}{n \brace k}_{f,\lambda}(H)\big(\ap\big)^{k}a^{k}=\sum_{k=0}^{n}{n \brace k}_{f,\lambda}(H)\sum_{j=0}^{k}S_{f,\lambda}(k,j,H)\big(\fa\big)_{j,\lambda} \label{25}\\
&=\sum_{j=0}^{n}\bigg(\sum_{k=j}^{n}{n \brace k}_{f,\lambda}(H)S_{f,\lambda}(k,j,H)\bigg)\big(\fa\big)_{j,\lambda}.\nonumber
\end{align}
Therefore, by \eqref{24} and \eqref{25}, we obtain the following theorem.
\begin{theorem}
For $n,j\ge 0$, we have
\begin{equation*}
\sum_{k=j}^{n}S_{f,\lambda}(n,k,H){k \brace j}_{f,\lambda}(H)=\delta_{n,j},\quad\textrm{and}\quad \sum_{k=j}^{n}{n \brace k}_{f,\lambda}(H)S_{f,\lambda}(k,j,H)=\delta_{n,j},
\end{equation*}
where $\delta_{n,j}$ is Kronecker delta.
\end{theorem}
Now, we observe \eqref{12} and \eqref{18} that
\begin{equation}\label{26}
\begin{aligned}
\big(\ap\big)^{n}a^{n}\big|m\big\rangle &=\sum_{k=0}^{n}S_{f,\lambda}(n,k,H)\big(\fa\big)_{k,\lambda}\bigg|m\bigg\rangle\\
&=\sum_{k=0}^{n}S_{f,\lambda}(n,k,H)\Big([m]_{f}(\alpha_{0})\Big)_{k,\lambda}\bigg|m\bigg\rangle,
\end{aligned}
\end{equation}
where
\begin{equation*}
\begin{aligned}
&\big([m]_{f}(\alpha_{0})\big)_{0,\lambda}=1,\\
&\big([m]_{f}(\alpha_{0})\big)_{n,\lambda}=[m]_{f}(\alpha_{0})\big([m]_{f}(\alpha_{0})-\lambda\big)\cdots \big([m]_{f}(\alpha_{0})-(n-1)\lambda\big),\ (n\ge 1).
\end{aligned}
\end{equation*}
On the other hand, by \eqref{13} and \eqref{14}, we get
\begin{equation}\label{27}
\begin{aligned}
\big(\ap\big)^{n}a^{n}\big|m\big\rangle &=[m]_{f}(\alpha_{0})[m-1]_{f}(\alpha_{0})\cdots [m-n+1]_{f}(\alpha_{0})\big|m\big\rangle \\
&=[m]_{f,n}(\alpha_{0})\big|m\big\rangle,
\end{aligned}
\end{equation}
where
\begin{equation*}
\begin{aligned}
&\big[m\big]_{f,k}(\alpha_{0})=[m]_{f}(\alpha_{0})[m-1]_{f}(\alpha_{0})[m-2]_{f}(\alpha_{0})\cdots [m-k+1]_{f}(\alpha_{0}),\ (k\ge 1) \\
&\textrm{and}\quad [m]_{f,0}(\alpha_{0})=1.
\end{aligned}	
\end{equation*}
Therefore, by \eqref{26} and \eqref{27}, we obtain the following theorem.
\begin{theorem}
For $n\ge 0$, we have
\begin{equation*}
[m]_{f,n}(\alpha_{0})=\sum_{k=0}^{m}S_{f,\lambda}(n,k,H)\Big([m]_{f}(\alpha_{0})\Big)_{n,\lambda}.
\end{equation*}
\end{theorem}
For $r\in\mathbb{N}\cup\{0\}$, the degenerate $r$-Stirling numbers of the first kind are defined by
\begin{equation*}
(x)_{n}=\sum_{k=0}^{n}S_{1,\lambda}(n+r,k+r)(x+r)_{n,\lambda},\quad (\mathrm{see}\ [7,8,11]).
\end{equation*}
Now, we define the {\it{degenerate generalized $r$-Stirling operators of the first kind}} by
\begin{align}
\big(N_{f}\big)_{n}=\big(\fa\big)_{n}&=\sum_{k=0}^{n}S_{f,\lambda}^{(r)}(n+r,k+r,H)\big(N_{f}+[r]_{f}(H)\big)_{k,\lambda}\label{28}\\
&=\sum_{k=0}^{n}S_{f,\lambda}^{(r)}(n+r,k+r,H)\big(\fa+[r]_{f}(H)\big)_{k,\lambda},\quad (n\ge 0). \nonumber	
\end{align}
Now, we observe from \eqref{9} that
\begin{align*}
\fa(\ap)^{r}&=\ap \big((\ap)^{r}a+(\ap)^{r-1}[r]_{f}(H)\big)  \\
&=(\ap)^{r}\big(\fa+[r]_{f}(H)\big),
\end{align*}
which implies that
\begin{align}
\big(\fa\big)_{n,\lambda}\big(\ap\big)^{r}&=\big(\fa\big)\big(\fa-\lambda\big)\cdots\big(\fa-(n-1)\lambda\big)\big(\ap\big)^{r} \label{29} \\
&=\big(\ap\big)^{r}\big(\fa+[r]_{f}(H)\big)_{n,\lambda}. \nonumber	
\end{align}
From \eqref{19}, \eqref{28} and \eqref{29}, we note that
\begin{align}
\big(\ap\big)^{r}(N_{f})_{n}\big|m\big\rangle
&=\big(\ap\big)^{r}\sum_{k=0}^{n}S_{f,\lambda}^{(r)}(n+r,k+r,H)\big(\fa+[r]_{f}(H)\big)_{k,\lambda}\big|m\big\rangle \label{30}\\	
&=\sum_{k=0}^{n}S_{f,\lambda}^{(r)}\big(n+r,k+r,f^{-r}(H)\big)\big(\ap\big)^{r}\big(\fa+[r]_{f}(H)\big)_{k,\lambda}\big|m\big\rangle \nonumber \\
&=\sum_{k=0}^{n}S_{f,\lambda}^{(r)}\big(n+r,k+r,f^{-r}(H)\big)\big(\fa\big)_{k,\lambda}\big(\ap\big)^{r}\big|m\big\rangle\nonumber
\end{align}
Therefore, by \eqref{30}, we obtain the following theorem,
\begin{theorem}
For $n,r\ge 0$, we have
\begin{equation*}
(\ap)^{r}(\fa)_{n}= \sum_{k=0}^{n}S_{f,\lambda}^{(r)}\big(n+r,k+r,f^{-r}(H)\big)\big(\fa\big)_{k,\lambda}\big(\ap\big)^{r}.
\end{equation*}
\end{theorem}

\section{Conclusion}
In this paper, we successfully extended the intersection of combinatorial calculus and quantum mechanics by investigating the degenerate generalized Stirling operators of the first kind within the framework of the generalized Heisenberg algebra $\text{GHA}_f$. While previous literature focused heavily on normal ordering expansions and Stirling operators of the second kind, this work provides the necessary inverse relations to complete the operational framework. Specifically, we established: An explicit factorization formula for the monomial operators $(\ap)^k a^k$. A fundamental triangular recurrence relation governed by the structural function $f(H)$ and the degenerate parameter $\lambda$. Strict orthogonality relations that cleanly map the transformations between the degenerate generalized Stirling operators of both kinds. A generalized shifting identity involving auxiliary operator indices.These results offer a robust mathematical framework for handling normal ordering problems in deformed quantum systems. Future work may explore the physical applications of these operational identities in specific quantum models, such as $q$-deformed oscillators, or investigate their implications for generalized coherent states.

\section{Conflict of interes}
 On behalf of all authors, the corresponding author states that there is no conflict of interest.

\end{document}